\documentclass[a4paper]{article}
\usepackage{latexsym,amsfonts,amscd,amssymb,enumerate,amstext,amsmath,dsfont,wasysym}
\usepackage[sc,tiny,center,compact]{titlesec}
\usepackage{amsthm}
\usepackage{graphicx}
\usepackage[all]{xy}
\usepackage{titletoc}
\titlecontents{section}[0pt]{\small\scshape}{\contentslabel{1.5em}}{}{\titlerule*[0.75pc]{.}\contentspage}[\textbullet\ ]

\titlecontents{subsection}[1.5em]{\small\scshape}{\contentslabel{2.3em}}{}{\titlerule*[0.75pc]{.}\contentspage}[\textbullet\ ]
\titlecontents{subsubsection}[3.8em]{\small\scshape}{\contentslabel{3.0em}}{}{\titlerule*[0.75pc]{.}\contentspage}[\textbullet\ ]

\author{M. Grime \\ University of Bristol\\ email: Matt.Grime@bris.ac.uk}
\title{Finite dimensional modules and perpendicular subcategories\footnote{2000 MSC 16G60 18A30 18Exx 18G25 20C05 20C20}}
\date{June 2007}
\newtheoremstyle{ordinary}{1ex}{0pt}{}{}{\scshape}{.}{\newline}{}
\theoremstyle{ordinary}
\newtheorem{thm}{Theorem}[section]
\newtheorem{defn}[thm]{Definition}
\newtheorem{lem}[thm]{Lemma}

\newtheorem{prop}[thm]{Proposition}

\newcommand{\tind}{\,\raisebox{0.6mm}{\rotatebox[origin=c]{90}{$\rightsquigarrow$}}}

\newcommand{\catt}{\mathcal{T}}

\newcommand{\cats}{\mathcal{S}}

\newcommand{\catp}{\mathcal{P}}
\renewcommand{\mod}{\mathrm{mod}(kG)}
\newcommand{\Mod}{\mathrm{Mod}(kG)}
\newcommand{\stmod}{\mathrm{stmod}_w(kG)}
\newcommand{\Stmod}{\mathrm{StMod}_w(kG)}

\begin{document}
\maketitle

\begin{abstract}
We explain how, under some hypotheses, one can construct a sequence of finite dimensional $kG$-modules that lie in certain prescribed additive subcategories, but whose direct limits do not. We use these to show that many of the triangulated quotients of $\Mod$ are not generated, as triangulated categories, by the corresponding quotient of $\mod$ considered as a full subcategory.

\end{abstract}

\section{Introduction}

Let $G$ be a finite group, and $k$ an field such that char$(k)$ divides $|G|$. The categories $\mod$ and $\Mod$ are Frobenius categories (see \cite{happel} for example, for an explanation), which implies that the quotients

\[ \mathrm{stmod}(kG):= \frac{\mod}{ \mathrm{f.g. projective }\ kG\mathrm{-modules}}\]
and
\[ \mathrm{StMod}(kG) \frac{\Mod}{\mathrm{projective }\ kG\mathrm{-modules}}\]
are triangulated categories. Whenever one has a triangulated category it is natural to ask if there is a smaller subcategory which generates it. Recall that if $\cats \subset \catt$ are triangulated categories, then $\cats$ generates $\catt$ if $(S,X)_\catt=0$ for all $S\in \cats$ implies that $X=0$. It is not too hard to show that $\mathrm{stmod}(kG)$ generates $\mathrm{StMod}(kG)$. Our aim is to show that in other triangulated quotients of $\Mod$ the finite dimensional objects do not often form a generating subcategory.  We will do this by producing a sequence of finite dimensional modules in $\mod$ that are zero in the quotient, but with direct limit in $\Mod$ that does not become zero.

%\section{Preliminary material}
%\subsection{Triangulated categories}
%We start by recalling some definitions.

%\begin{defn}[Compact objects]
%Let $\catt$  be a triangulated category, an object $c$ in $\catt$ is said to be compact (or small) if the canonical map

%\[ \coprod(c,X_\alpha) _\catt\to (c,\coprod X_\alpha)_\catt\]
%is an isomorphism. Equivalently, $c$ is compact if every map into a coproduct factors through finitely many terms. 
%\end{defn}
%\begin{defn}[Generation of triangulated categories]
%Let $\cats$ be a class of objects  in $\catt$. To $\cats$ we can associate a thick, localizing subcategory, $\cats^\perp$ called the orthogonal category:

%\[\cats^\perp:= \{ X : (\cats,X)_\catt = 0\}\]
%We say $\cats$ generates $\catt$ if $\cats^\perp=0$. 
%\end{defn}

%\begin{defn}[Compact generation]
%A triangulated category $\catt$ is said to be compactly generated if there is a \emph{set} of compact objects that generate $\catt$.
%\end{defn}

\section{Modular representation theory and triangulated quotients}
We continue with the assumption that $k$ is a field, and char$(k)$ divides $|G|$.  We assume that the reader is familiar with the content of, say, Alperin's book \cite{alperin}.
\begin{defn}[Relatively projectivity]
Let $w$ be a finite dimensional $kG$-module. Let $\catp(w)$ denote the smallest additive subcategory of $\Mod$ that contains $w$ and is closed under tensor with an arbitrary module and arbitrary direct sums and summands.
\end{defn}
The class $\catp(w)$ is sufficient to allow a relative cohomology theory, and a triangulated quotient of $\Mod$. Objects in $\catp(w)$ are called $w$-projective.

\begin{thm} Let $\Delta$ be the class of short exact sequences in $\Mod$ that split when tensored with $w$. Then $\Delta$ is an exact structure on $\Mod$, and the class of objects $\catp(w)$ constitute the projective and injective objects with respect to that structure. Moreover, there are enough pro/injective objects, and we can define triangulated quotients

\[\Stmod:= \frac{\Mod}{\catp(w)} \ \ \ \   \stmod:=\frac{\mod}{\catp\cap\mod}\]
\end{thm}
\begin{proof}
See, e.g. \cite{virtual}.
\end{proof}

If one picks a subgroup $H<G$, and sets $w=\mathrm{Ind}_H^G(k)$, then one obtains the usual definition of $H$-projective. The ordinary stable category can be recovered by choosing $w=kG$.

\subsection{Twisting $kH$-modules}\label{twist}

We continue with the assumption that $k$ is a field of characteristic $p$, and further suppose that $q$ is a power of $p$. Let $H$ be a group,
\[\xymatrix{ 0\ar[r]& X \ar[r]^{d_1} & Y \ar[r]^{d_2} & Z \ar[r] & 0 }\]
a short exact sequence of $kH$-modules, and let $G=H\times C_q$. We wish to use this short exact sequence to define a $kG$-module, $(X,Y,Z)\tind_H^G$. The reader should think of $\tind$ as meaning \emph{twisted induction}\footnote{I am indebted to Jeremy Rickard for suggesting this construction to me} .  As a vector space sum $(X,Y,Z)\tind_H^G$ will be given by
\[ \underbrace{X +\cdots + X}_{q-1} + Y + \underbrace{Z+\cdots + Z}_{q-1}\]
and the $H$-action will be the obvious one in each summand. Thus it remains to describe the $C_q$-action. Let $C_q$ be generated by $u$. Then $u-1$ acts by shifting summands in the following manner:
\[\xymatrix{ X \ar[r]^1 & {\cdots} \ar[r]^1 & X \ar[r]^{d_1} & Y \ar[r]^{d_2} & Z \ar[r]^1 & {\cdots} \ar[r]^1 & Z\ar[r]^0 & 0}\]
Notice that $(u-1)^q=u^q-1=0$, since applying $u-1$ $q$ times to any of the summands of $(X,Y,Z)\tind_H^G$ will mean applying $d_2d_1$, or $0$, at some point.
\begin{prop} Let $X,Y,Z$ and $(X,Y,Z)\tind_H^G$ be as above, then $(X,Y,Z)\tind_H^G$ is $H$-projective if and only if the map $X\to Y$ splits.
\end{prop}
\begin{proof}
 Consider $X\otimes k$ as a $kH\times C_q$ module with the diagonal action. The module $(X,Y,Z)\tind_H^G$ is $H$-projective if and only if it has the lifting property with respect to $H$-split short exact sequences, thus consider the $H$-split surjection
\[ \pi:\mathrm{Ind}_H^G( \mathrm{Res}_H^G(X\otimes k)) \twoheadrightarrow  X\otimes k\]
There is a map from $(X,Y,Z)\tind_H^G$ to $X\otimes k$ given by projection into the first copy of $X$. We will show that this map factors through $\pi$ if and only if $d_1: X\to Y$ is split.

Suppose that $\theta$ is such that $(X,Y,Z)\tind_H^G\to X\otimes k$ factors as $\pi\theta$.  We will use the vector space decomposition of $(X,Y,Z)\tind_H^G$ as above, and we can consider $\mathrm{Ind}_H^G(X)$ as a vector space sum
\[ \underbrace{X + \ldots +X}_q\]
in the usual manner: the summands are indexed by cosets of $H$, i.e. powers of $u$. Let us write $\theta$ as a block  matrix with respect to these vector space sum

\[\theta = \left( \begin{array}{cccc} \theta_{1,1} & \theta_{1,2} & \cdots & \theta_{1,2q-1} \\ \vdots & \vdots & \vdots &\vdots \\ \theta_{q,1} & \theta_{q,2}& \cdots & \theta_{q,2q-1}   \end{array} \right)\]
with each $\theta_{r,s}$ a $kH$-equivariant map. Similarly $\pi$ is given by 

\[ (\underbrace{1,\ldots,1}_q)\]
and the projection from $(X,Y,Z)\tind_H^G$ to $X\otimes k$ is
\[(1,\underbrace{0,\ldots,0}_{2q-2})\]
From the factorization of the projection as $\pi\theta$ one obtains

\[ ( \sum_i \theta_{i,1} , \sum_i \theta_{i,2},\ldots, \sum_i \theta_{i,2q-1}) = (1,0,\ldots,0).\]
We also know that $h\theta=\theta h$.  The reader is encouraged to work out the case of $q=2$ by hand, and to write down the matrices for larger $q$. When they have done so they will notice that one has the extra relations (indices are to be read mod $q$)

\begin{eqnarray}\theta_{r,s} &=&\theta_{r+1,s}+\theta_{r+1,s+1}  \ \ \  1\leq r \leq q-1\label{relone} \\
 \theta_{q,q-1}  &=&\theta_{1,q-1} + \theta_{1,q}d_1 \label{reltwo} \end{eqnarray}
It follows from (\ref{relone}) by induction on $k$ that

\[ \theta_{r,s} = \sum_{i=0}^{k} \binom{k}{i}\theta_{r-k+i,s-k}\]
for all $1\leq s\leq q-1$. We find it easiest not to insert limits in the sums for what follows. Recall that we define $\binom{n}{m}$ to be zero if $m$ is not between $0$ and $n$. Thus, using the relations we generated, (\ref{reltwo}), and showing a healthy disregard for indices it follows that

\begin{eqnarray*}\theta_{1,q}d_1 &=& \theta_{q,q-1} - \theta_{1,q-1}\\
& = &\sum_{i}(-1)^i\binom{q-1}{i}\theta_{p-(p-1)+i,1} - \sum_{i} (-1)^i\binom{q-1}{i}\theta_{1-(q-1)+i,1} \\
& =&  \sum_{i}(-1)^i\binom{q-1}{i}\theta_{1+i,1} - \sum_{i} (-1)^i\binom{q-1}{i}\theta_{2-q+i,1}\\
 &=& \sum_{i}(-1)^i\binom{q-1}{i}\theta_{1+i,1} - \sum_{i} (-1)^i\binom{q-1}{i}\theta_{2+i,1}\\
 &=&\sum_{i}(-1)^i\left(\binom{q-1}{i}+\binom{q-1}{i-1}\right)\theta_{1+i,1}\\
 &=&\sum_{i} (-1)^i\binom{q-1}{i}\theta_{1+i,1} = \sum_{i} ((-1)^i)^2\theta_{1+i,1}\end{eqnarray*}
and thus (recalling that indices are mod $q$)
\[ \theta_{1,q}d = \sum_i \theta_{i,1} = 1_X \]
which completes the proof that $d_1$ splits.
\end{proof}

\section{The inclusion of $\stmod$ in $\Stmod$}

 In this section we will argue that under some reasonable assumptions on $G$, and $w$, we may show that $\stmod$ does not generate $\Stmod$. The tactic is to write some non-$w$ projective module as a direct limit of $w$-projective modules. In fact, we shall show something slightly stronger: the direct limit will not be $\mathrm{vtx}(w)$-projective.  First, we will need a way to show a module is not $w$-projective.
 
 \begin{lem}
 Let $X$ be a $w$-projective $kG$-module, then $X$ is projective with respect to any vertex of $w$. 
 \end{lem}
 \begin{proof}
 It suffices to consider the case $X\cong w\otimes Y$. Let $Q$ be a vertex of $w$ and let $v$ be a source. Then 
 
 \[ w\otimes Y \ | \ \mathrm{Ind}_Q^G(v) \otimes Y \cong \mathrm{Ind}( v\otimes \mathrm{Res}^G_Q(Y))\]
 and we see $X$ is $Q$-projective.
\end{proof} 
 Now we show that it suffices to pass to the Sylow-$p$ subgroup of $G$. 
 
 \begin{prop} Suppose that $P$ is a Sylow-$p$ subgroup of $G$ and let $v$ the restriction of $w$ to $kP$.  Suppose that $M= \varinjlim m_\alpha$ is a filtered colimit in $\mathrm{Mod}(kP)$  where each $m_\alpha$ is finite dimensional and  $v$-projective and $M$ is not projective with respect to $\mathrm{vtx}(w)$, then  $\mathrm{Ind}_P^G(M) = \varinjlim \mathrm{Ind}_P^G(m_\alpha)$  is a non-$w$-projective $kG$-module that is the direct limit of finite dimensional $w$-projectives.
 \end{prop}
 \begin{proof}
 This is reasonably clear by the last lemma.
 \end{proof}
 Thus we may suppose that $G$ is a $p$-group. The most natural statement (i.e.~the one with fewest hypotheses) is when $w=\mathrm{Ind}_H^G(k)$.

 \begin{thm}\label{pgroup} Let $H$ be a $p$-group with non-finite representation type, and let $G=H\times C_q$. Set $w=\mathrm{Ind}_H^{G}(k)$, then $\Stmod$ is not generated by $\stmod$.
\end{thm}
\begin{proof}

 The hypothesis on $H$ ensures that there is an indecomposable countable dimensional $kH$-module $M$. Suppose that we write $M$ as the direct limit of a sequence of finite dimensional modules

\[ \varinjlim_{n \in \mathbb{N}} m_n\]
Let $\iota_n$ denote the inclusion of $m_n$ into $m_{n+1}$ and consider the (non-split) short exact sequence

\[\xymatrix{ 0\ar[r] & {\coprod} m_n \ar[r]^{1-\iota_n} & {\coprod} m_n \ar[r] & M \ar[r] & 0 }\]
which is the direct limit of the split short exact sequences

\[  0 \to  \mathop{\coprod}_{n=1}^{N} m_n  \to \mathop{\coprod}_{n=1}^{N+1} m_n \to m _{N+1} \to 0 \]
Construct the module $(\coprod m_n,\coprod m_n,M)\tind_H^G$ as in subsection \ref{twist}. This is \emph{not} $H$-projective as the map $\coprod m_n \to M$ does not split. However, $(\coprod m_n,\coprod m_n,M)\tind_H^G$ is the direct limit of the modules 

\[  (\mathop{\coprod}_{n=1}^{N} m_n  ,\mathop{\coprod}_{n=1}^{N+1} m_n , m _{N+1} )\tind_H^G \]
each of which is $H$-projective, since  $\mathop{\coprod}_{n=1}^{N+1} m_n \to m _{N+1}$ is split. Now, any map from a finite dimensional $kG$-module to $(\coprod m_n,\coprod m_n,M)\tind_H^G$ factors through a finite dimensional submodule, and thus through some $H$-projective submodule. Hence $(\coprod m_n,\coprod m_n,M)\tind_H^G$ is orthogonal to the set of finite dimensional modules.
\end{proof}
All that remains is to extend this to the case when $w$ is not a trivial source module.
\begin{thm}\label{intermed} Let $G$, $H$ and $M$ be as in \ref{pgroup}, and suppose that $w$ is a $kG$-module with $\mathrm{vtx}(w) \leq_G H$. Suppose that $\mathrm{Res}_H^G(w)\otimes M$ is not pure projective. Then $\stmod$ does not generate $\Stmod$.
\end{thm}

\begin{proof}
 Let $v=\mathrm{Res}_H^G(w)$. The hypotheses imply that 
\[\xymatrix{ 0\ar[r] & {\coprod} v\otimes m_n \ar[r]^{1-\iota_n} & {\coprod} v \otimes m_n \ar[r] & v\otimes  M \ar[r] & 0 }\]
is non-split, and hence, $(\mathop{\coprod} v\otimes m_n  ,\mathop{\coprod} v \otimes m_n , v \otimes M )\tind_H^G$ is not $w$-projective. This module is the direct limit of the modules

\[(\mathop{\coprod}_{n=1}^{N}v \otimes m_n  ,\mathop{\coprod}_{n=1}^{N+1}v\otimes m_n , v\otimes m _{N+1} )\tind_H^G. \]
Thus we need to show that these are $w$-projective for each $n$. We know that any such module is $H$-projective, which means it is a summand of

\[ \mathrm{Ind}^G_H(\mathrm{Res}_H^G((\mathop{\coprod}_{n=1}^{N} v\otimes m_n  ,\mathop{\coprod}_{n=1}^{N+1} v \otimes m_n , v \otimes m _{N+1} )\tind_H^G)) \]
but this is nothing more than a direct sum of copies of modules of the form

\[  \mathrm{Ind}_H^G(v \otimes m_n) \cong w \otimes \mathrm{Ind}_H^G(m_n)\]
and thus is $w$-projective as we were required to show.

\end{proof}

We will end by collecting these theorems into one statement.
\begin{thm} Let $G$ be a finite group and $w$ a $kG$-module with vertex (conjugate to) $Q$. Assume that $G$ and $w$ satisfy the following conditions:
 
 \begin{itemize}
 \item  if $P$ is a Sylow-$p$ subgroup of $G$ then $P$ is isomorphic to $P'\times C_q$ for some $P'$;
 \item in such a decomposition $ Q \leq_GP'$;
  \item there is a $kP'$ module $M$ such that neither $M$ nor $M\otimes \mathrm{Res}(w)$ are pure projective (in particular $kQ$ and thus  $kP'$ cannot have finite representation type);
 \end{itemize}
 then $\stmod$ does not generate $\Stmod$ as a triangulated category.

\end{thm}

 \begin{small}

\end{small}

\end{document}